\documentclass[lettersize,journal]{IEEEtran}
\usepackage{amsmath,amsfonts}
\usepackage{algpseudocode}
\usepackage{algorithm}
\usepackage{array}
\usepackage[caption=false,font=normalsize,labelfont=sf,textfont=sf]{subfig}
\usepackage{multirow}
\usepackage{textcomp}
\usepackage{stfloats}
\usepackage{url}
\usepackage{verbatim}
\usepackage{graphicx}
\usepackage{cite}
\usepackage{hyperref}
\usepackage[table]{xcolor}
\usepackage{enumitem}
\begin{document}

\title{Regularized Benders Decomposition for High Performance Capacity Expansion Models}

\author{Filippo Pecci\thanks{Filippo Pecci is with the Andlinger Center for Energy and the Environment, Princeton University, Princeton, NJ, United States of America (e-mail: \text{filippopecci@princeton.edu}).}, Jesse D. Jenkins\thanks{Jesse D. Jenkins is with the Andlinger Center for Energy and the Environment and the Department of Mechanical and Aerospace Engineering, Princeton University, Princeton, NJ, United States of America (e-mail: \text{jessejenkins@princeton.edu}).

This work has been accepted for publication in IEEE Transactions on Power Systems. © 2025 IEEE. Personal use of this material is permitted. Permission from IEEE must be obtained for all other uses, in any current or future media, including reprinting/republishing this material for advertising or promotional purposes, creating new collective works, for resale or redistribution to servers or lists, or reuse of any copyrighted component of this work in other works.
}}
\markboth{IEEE Transactions on Power Systems}{Pecci and Jenkins \MakeLowercase{\textit{et al.}}: Regularized Benders Decomposition for High Performance CEMs}


\maketitle

\begin{abstract}
    We consider electricity capacity expansion models, which optimize investment and retirement decisions by minimizing both investment and operation costs. In order to provide credible support for planning and policy decisions, these models need to include detailed operations and time-coupling constraints, {consider multiple possible realizations of weather-related parameters and demand data,} and allow modeling of discrete investment and retirement decisions. Such requirements result in large-scale mixed-integer optimization problems that are intractable with off-the-shelf solvers. Hence, practical solution approaches often rely on carefully designed abstraction techniques to find the best compromise between reduced {computational burden} and model accuracy. Benders decomposition offers scalable approaches to leverage distributed computing resources and enable models with both high resolution and computational performance. 
    In this study, we implement {a tailored Benders decomposition method for large-scale capacity expansion models with multiple planning periods, stochastic operational scenarios, time-coupling policy constraints, and multi-day energy storage and reservoir hydro resources. Using multiple case studies, we also} evaluate several level-set regularization schemes to {accelerate convergence. We} find that a regularization scheme that selects planning decisions in the interior of the feasible set shows superior performance compared to previously published methods, enabling high-resolution, mixed-integer planning problems with unprecedented computational performance.
\end{abstract}


\section{Introduction}
As the urgency to mitigate greenhouse gas emissions grows, capacity expansion models (CEMs) play a key role in providing decision support for a sustainable and equitable transition in the electricity sector, helping to avoid misallocation of investments, understand the role of emerging technologies and the impact of possible policy interventions, and develop robust strategies to transition to net-zero emissions electricity systems~\cite{sepulveda_design_2021,victoria_speed_2022,ricks_minimizing_2023}. 
Ideally, these models should co-optimize generation and transmission investment and retirement decisions across multiple planning periods, and operational dispatch decisions at hourly or finer resolution under different scenarios, with detailed operating constraints (e.g., unit commitment, energy storage), and high geospatial resolution~\cite{levin_energy_2023}. 
These requirements result in least-cost optimization problems that combine integer and continuous decision variables subject to a set of linear constraints representing engineering operation and environmental, political, and economic requirements. Integer variables are used to model discrete investment and retirement decisions (and associated economies of scale) across multiple planning periods, which are optimized to minimize capital and operating costs. To estimate operating costs, CEMs include a representation of the underlying electricity system, where continuous variables model dispatch and storage decisions, and linear constraints are used to model physical and operational limits of the considered technologies. The electricity network is modeled as a graph, whose nodes represent different geographical zones. Power flows between zones are represented by a transport model with losses or a linearized optimal power flow model~\cite{neumann_assessments_2022}. 
{To improve solution robustness with respect to uncertain future weather conditions, each planning period should correspond to multiple realizations of weather-related parameters and demand data, modeled with hourly resolution (i.e., multiple stochastic scenarios).}
Therefore, full resolution CEMs result in large-scale mixed-integer linear programs (MILPs) with tens to hundreds of millions of variables and constraints, pushing even the best commercial solvers to their computational limits. 
Due to these computational constraints, existing models are heavily simplified by: sampling representative time periods or ignoring sequential operations entirely (e.g., ``time slices''); aggregating regions into larger geographic zones; and/or ignoring or relaxing key operational constraints. These abstractions can ensure CEMs are computationally tractable but come at the cost of significantly reduced accuracy that impacts their ability to provide credible decision support. As examples, see~\cite{palmintier_heterogeneous_2014} on the impacts of operational simplification, the discussion in~\cite{frysztacki_strong_2021} on the effects of spatial aggregation, and the studies~\cite{poncelet_impact_2016,mallapragada_impact_2018,jacobson_quantifying_2024} on the potential biases introduced by temporal aggregation. {In particular, the numerical results in~\cite{jacobson_quantifying_2024} indicate that the use of few representative days can significantly affect both renewables siting decisions as well key operational outcomes like demand curtailment and CO$_2$ emissions. Moreover,} modelling of multi-day energy storage (MDS) resources is becoming central in power system planning, especially when considering the interactions between electricity systems and other energy carriers, or key industrial sectors~\cite{levin_energy_2023}. However, accurate representation of MDS requires a higher temporal resolution than what is usually considered in capacity expansion models.  

As an alternative to abstractions, several studies have focused on developing decomposition methods for solving large optimization problems arising in the framework of energy system planning. A common approach is to implement Benders decomposition to separate planning decisions from operational dispatch decisions. At every iteration, this cutting plane algorithm requires solving two optimization problems: (i) a planning problem to select investment and retirement decisions; (ii) a set of operational sub-problems optimizing dispatch decisions over each planning period to generate cuts (linear inequalities) to be added to the planning problem. In this standard setting, an operational sub-problem represents a full year of system operation, modeled with hourly resolution, coupled in time by both operational and policy constraints. To reduce computational cost, some studies do not model integer investment and retirement decisions~\cite{goke_stabilized_2024}, and simplify the operational model by not considering storage resources~\cite{munoz_new_2016}, or
ignoring ramping constraints and unit commitment~\cite{munoz_new_2016,zhang_integrated_2023,goke_stabilized_2024}. Similarly to these previous works, the study in~\cite{zhang_integrated_2023} considers decomposition schemes that separate planning problem and full-year operational sub-problem. To further reduce computational burden, they propose to avoid solving a full operational sub-problem at every iteration, using inexact cuts based on the approximate solutions of earlier exact iterations.  Such approach results in significant computational savings, but it requires the operational sub-problems from different planning periods to have the same constraint matrices, with variations allowed only in the right-hand side of the inequalities. 
For this reason, the capacity expansion models formulated in~\cite{zhang_integrated_2023} consider identical renewable availability profiles across planning periods, as well as the same demand hourly patterns. They represent varying demand across planning periods by applying a scaling factor to the demand hourly profile. This limits their application when considering different weather scenarios (i.e., renewables availability patterns). Further investigations are needed to assess the efficacy of these inexact cuts when variable parameters include hourly demand and renewable availability time series, like the general form in Problem~\eqref{eq:probform}.
The works in \cite{Lara2018,li_mixed-integer_2022} included detailed operational and time coupling constraints, but then formulated the operational sub-problem only for a selection of representative sub-periods from each year and do not model multi-day energy storage or hydropower reservoir resources (which introduce time coupling across operational sub-periods). 
{In comparison, \cite{soares_integrated_2022} represents hydro-reservoirs with monthly temporal resolution, ignoring their chronological hourly operation. The rest of the electricity system operation is represented with hourly resolution only for 3 typical days per month. While these approximations are acceptable in systems dominated by hydropower and thermal generation, they might introduce significant biases when considering large shares of intermittent weather-dependent renewable energy sources (e.g., wind and solar power), whose hourly generation significantly varies during the day and between days of the week - see~\cite{poncelet_impact_2016}, \cite{mallapragada_impact_2018}, and \cite{jacobson_quantifying_2024}.}

{Compared to existing literature, we consider CEMs with multiple planning periods and operational scenarios, hourly temporal resolution, unit commitment and ramping constraints, integer planning decisions, time-coupling policy constraints, and multi-day energy storage resources. We develop a tailored Benders decomposition framework for this class of CEMs, able to decompose full operational years with hourly resolution into sub-problems that can be optimized in parallel. In a recent study,} \cite{jacobson_computationally_2024} presented a model reformulation and solution algorithm for single-period CEMs (and a single operational scenario) that allows us to decompose the full operational year into shorter sub-periods. {The study focused on CEMs with time-coupling policy constraints but excluded multi-day energy storage resources, introducing an equivalent reformulation based on policy budgeting variables to link different sub-periods, where the budgeting variables are optimized by the planning problem.} By generating multiple cuts per iteration (one for each operational sub-period), the method accelerates convergence relative to alternative Benders formulations that produce fewer cuts per iteration (one per each full operational year)~\cite{munoz_new_2016,Lara2018,li_mixed-integer_2022,soares_integrated_2022, zhang_integrated_2023,goke_stabilized_2024} {- see Figure 5 in \cite{jacobson_computationally_2024}}. Moreover, by taking advantage of distributed computing resources to solve the operational sub-problems in parallel, this algorithm significantly reduces the computational cost associated with each iteration compared to other literature. {However, the formulation in~\cite{jacobson_computationally_2024} did not consider multi-period CEMs, and it can not be applied when long duration (e.g., multi-day) energy storage or reservoir hydropower resources are present, because they introduce additional time-coupling between sub-periods, and thus prevent parallelization. In this study, we extend the formulation in~\cite{jacobson_computationally_2024} by considering as planning decisions both policy budgeting variables and initial and final storage levels for long duration energy storage and hydropower reservoir resources. This allows us to decompose a full-year operational sub-problem into smaller parallelized sub-problems, even in the presence of multi-day storage constraints.

Furthermore, as shown in Section~\ref{sec:numexp}, the Benders decomposition algorithm in~\cite{jacobson_computationally_2024} can suffer from slow convergence due to oscillations between extreme planning decisions.} This is a known issue of cutting plane methods, exacerbated by the large-number of planning decision variables included in multi-period CEMs. In continuous optimization, regularization techniques like {proximal-bundle~\cite{ruszczynski_regularized_1986}} and trust-region ~\cite{linderoth_decomposition_2003,goke_stabilized_2024} {methods} mitigate this oscillating behavior by controlling the step size at each iteration. The presence of discrete investment and retirement decisions complicates the use of proximal penalty terms and trust-regions, since integer feasible solutions may be far apart in the decision space. 
An alternative approach is offered by level-set methods\cite{lemarechal_new_1995}, where each iteration solves a regularization problem selecting a sub-optimal feasible solution of the planning problem, with a bound on its level of sub-optimality controlled by a tunable parameter. {Common implementations select the sub-optimal feasible solution by minimizing the distance to the current best solution or to the minimizer of the non-regularized planning problem~\cite{lemarechal_new_1995}. Alternatives include interior-point level-set regularization~\cite{gondzio_new_2013,zhang_integrated_2023}, which aims to select a sub-optimal planning solution that belongs to the interior of the feasible set of the planning problem, with the level-set constraint used to control the degree of sub-optimality.

Compared to proximal-bundle and trust-region methods, level-set methods depend on a single parameter, which does not need to be dynamically updated at each iteration. 
Moreover, the experimental study~\cite{zverovich_computational_2012} reported that a level set regularization based on the $\ell_2$-norm outperformed both a trust-region scheme and a proximal-bundle method on a library of benchmarking linear programs. In the context of capacity expansion models, the work in~\cite{goke_stabilized_2024} has reached different conclusions, reporting a better computational performance for a trust-region regularization approach on the single system they examined. Both~\cite{goke_stabilized_2024} and~\cite{zverovich_computational_2012}
did not consider the interior-point level-set regularization, which outperforms all other strategies in the numerical experiments reported in Section~\ref{sec:numexp}.

Previous studies applying regularization to solve capacity expansion models (e.g., ~\cite{soares_integrated_2022,zhang_integrated_2023,goke_stabilized_2024}) considered a full-year operational sub-problem (rather than a set of smaller parallelized operational sub-problems), and the planning decision variables subject to regularization corresponded only to capacity investment and retirement decisions. In our framework, the set of planning decisions includes also continuous variables coupling sub-periods, namely budgeting variables for policy constraints, and initial and final storage levels for multi-day energy storage and hydropower reservoirs. Furthermore, level-set methods have usually been applied to continuous optimization problems. In particular, the interior-point level-set regularization loses meaning in a discrete feasible set that does not have an interior. Hence, we propose a two-stage algorithm to benefit from interior-point level-set regularization in presence of {mixed-integer} decision variables.  {In our approach, we initially implement an interior-point regularized Benders decomposition algorithm to solve a relaxation of the CEM with continuous planning decisions, enforcing integrality constraints only when the continuous relaxation has reached convergence. At this stage, we apply interior-point regularization only to continuous planning decisions, while discrete decisions are set by the solution of the planning problem.} The combination of interior-point level-set method and our Benders decomposition algorithm results in a capacity expansion model with high computational performance and level of detail (e.g., hourly operations, multi-day energy storage/reservoir hydro, discrete investment and retirement decisions). 

All methods are numerically evaluated using 4 systems with different sizes, based on the continental United States. The largest case study has 26 zones, over a thousand resources, three planning periods, and hourly temporal resolution for 52 weeks in each planning stage (26,208 total operational periods), which result in a CEM with 111 million variables (including 4,650 integer variables) and 197 million constraints. Using these test systems, we investigate different regularization techniques and compare them with solution algorithms previously proposed by \cite{Lara2018,li_mixed-integer_2022,jacobson_computationally_2024}. To further explore the generalizability of our methods for this class of problems, we also consider a case study from~\cite{deng_harmonized_2023}, including a single-period planning model for Brazil with 8 different weather-years (e.g., multiple scenarios for wind, solar, hydropower, and demand time series representing uncertainty about future weather conditions).
In summary, this manuscript moves beyond previous literature by introducing the following advances:
\begin{itemize}
    \item  {Extend the parallelized Benders decomposition method presented in \cite{jacobson_computationally_2024} by proposing an equivalent reformulation for multi-day energy storage and reservoir hydropower resources that allows to optimize each operational sub-period in parallel.}
    \item Combine our novel decomposed formulation with level-set regularization methods to accelerate the Benders decomposition algorithm,  making it computationally feasible to extend the method to much larger planning problems than previous studies, including multi-period planning models as well as single-period models with several weather and demand scenarios (e.g. stochastic operational scenarios) with hourly resolution.
    \item Propose a two-stage approach for solving CEMs with mixed-integer planning decisions using interior-point level-set regularization to generate valid Benders cuts and warm-start the mixed-integer planning problem while substantially reducing the number of iterations of the final MILP problem. 
    \item {Evaluate different regularization and decomposition schemes using case studies with varying sizes and levels of complexity, identifying the most effective technique for this class of problems.}
\end{itemize}}
\section{Multi-period capacity expansion models with stochastic operational scenarios}
\label{sec:multi_cems}
This study considers multi-period power system capacity expansion models, where investment and retirement decisions are optimized over multiple planning periods, while also optimizing operational decisions with hourly resolution. { In our formulation, each planning period can correspond to multiple weather and demand scenarios (i.e., multiple stochastic scenarios). Let $S$ be the index set of stochastic scenarios, and let $\beta_s>0$ be the probability assigned scenario $s \in S$. Each planning period $p \in P$ is comprised of sub-periods $W_p = \cup_{s\in S}W_{p,s}$, where $W_{p,s}$ is the subset of sub-periods for operational scenario $s \in S$. Denote by $\sigma(w)$ the index of the scenario that includes operational sub-period $w$, i.e. such that $w \in W_{p,\sigma(w)}$. Here, a sub-period represents an operational week (i.e., $168$ hours). In this way, a planning period consists of $|S| \times 52$ sub-periods ($|S| \times 8736$ hours).} CEMs are formulated as mixed-integer linear programs (MILPs) with the objective of minimizing total system cost. Integer variables represent generation and storage capacity investment and retirement as well as transmission expansion decisions. We denote the vector of investment and retirement decision variables as $y_p$ for each period $p \in P$.
Generator dispatch decisions, storage levels, power flows, and unit commitment decisions are included in the vector of continuous variables $x_w \in \mathbb{R}^{n_x}$ for each operational sub-period $w \in W_p$ (e.g., week). {For illustrative purposes, Figure~\ref{fig:decision_tree} shows an example of the considered multi-period planning decision process, which includes $3$ planning periods (years), $5$ operational sub-periods (weeks), and $2$ stochastic scenarios.}
{
\begin{figure}
    \centering
    \includegraphics[width=0.30\textwidth]{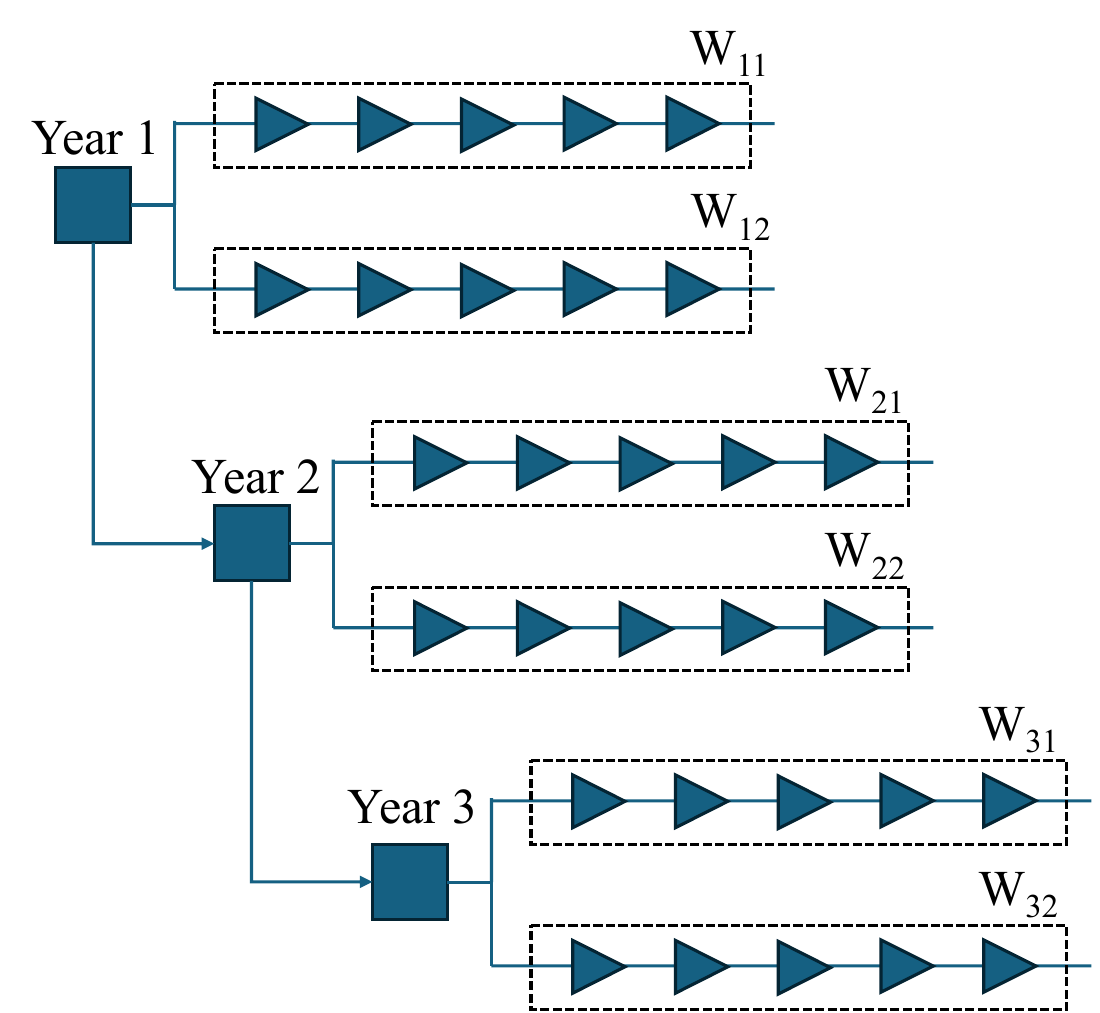}
    \caption{Illustrative example of the planning decision process with 3 planning years, each corresponding to 5 operational weeks and 2 stochastic scenarios. The operational weeks consist of 168 hours.}
    \label{fig:decision_tree}
\end{figure}
}
{We develop a Benders decomposition algorithm to decompose full operational years into shorter sub-periods, which can be optimized in parallel. Our framework is presented in Sections~\ref{sec:benders} and~\ref{sec:benders_reg}, and it applies to all CEMs that are represented by the following mixed-integer linear optimization problem:
\begin{subequations}
    \label{eq:probform}
    \begin{alignat}{3}
        &\text{min}&\; \; &\sum_{p \in P}\Big(f_p^Ty_p+\sum_{w \in W_p} \beta_{\sigma(w)} c_w^T x_w\Big) \label{eq:obj}\\
        &\text{s.t.}& & A_wx_w + B_wy_p \leq b_w, \quad \forall w \in W_p, \; p \in P \label{eq:opcons}\\
        &&& Q_w x_w \leq z_w, \quad \forall w \in W_p, \; p \in P \label{eq:auxcons1}\\
        &&& C_w z_w + D_wy_p \leq d_w, \quad \forall w \in W_p, \; p \in P \label{eq:auxcons2}\\
        &&& {\sum_{w \in W_{p,s}}E_wz_w = e_{p,s}, \quad \forall s \in S,\;\forall p \in P}\label{eq:couplingcons}\\
        &&& \sum_{p \in P}R_p y_p \leq r \label{eq:multiperiodcons} \\
        &&& x_w\geq 0,\quad \forall w \in W_p, \; \forall p \in P,\\
        &&& z_w \geq 0,\quad \forall w \in W_p, \; \forall p \in P,\\
        &&& y_p \geq 0, \quad \forall p \in P\\
        &&& y_{p,i} \in \mathbb{Z}, \quad \forall i \in I, \; \forall p \in P.
    \end{alignat}
\end{subequations}}
The objective function~\eqref{eq:obj} defines total system cost as the sum of fixed costs and variable costs, denoted here by vectors $f_p$ and $c_w$, respectively. Fixed costs include fixed operation and maintenance costs, as well as the cost of investment in new generation and transmission capacity. Variable costs consider fuel prices, variable operation and maintenance costs, and eventual penalties for non-served energy or other constraint violations. {The block diagonal} constraints~\eqref{eq:opcons} group operational constraints that apply separately to each sub-period $w \in W_p$, for all $p \in P$. These include energy balances in demand zones, capacity constraints for all power generation and storage resources, and transmission constraints. Finally, \eqref{eq:auxcons1}, \eqref{eq:auxcons2} and \eqref{eq:couplingcons} are used to represent constraints that couple different sub-periods, for example multi-day energy storage and policy constraints. {Note that, when multiple operational scenarios are considered, constraints~\eqref{eq:couplingcons} only link sub-periods from the same operational scenario.} {Constraints~\eqref{eq:multiperiodcons} couple planning decisions made over different periods, for example by enforcing that existing capacity in period $p$ is equal to the final capacity in period $p-1$. They also include constraints on minimum capacity retirements that must be made by the end of a given planning period.}

In the following, we discuss the most relevant properties of the class of capacity expansion models considered in this study.
In power system modeling, power flow is represented using energy networks, where different buses are connected via transmission lines. Because CEMs representing energy networks with thousands of nodes (i.e. buses) and transmission lines are computationally intractable without substantially reducing temporal resolution (e.g.~\cite{shawhan_costs_2019}), CEMs more commonly aggregate nodes within the same geographical area that have similar demand and climate conditions and between which power flows are rarely constrained. Similarly, we cluster resources within each aggregate node (or zone) based on their location, technology type, cost of connection to the grid, and operational properties. Therefore, power systems are represented by graphs whose nodes correspond to zones with aggregated demand and several aggregated generation and storage resource clusters, and edges model the interzonal transmission constraints, resulting from the aggregation of transmission lines. As observed in \cite{neumann_assessments_2022}, the inclusion of power flow equations, e.g. Kirchhoff’s Voltage Law (KVL), has a greater impact on models with lower level of aggregation, compared to highly aggregated models. Since the system considered in Section~\ref{sec:numexp} is a 26-zone model representing the continental United States, the case study does not consider KVLs when representing transmission between aggregated zones. We instead employ a lossy transport model that includes linearized power losses across interzonal transmission lines. However, the same methods and formulations discussed in this paper can be applied to models with linearized power flow equations (e.g., KVLs). Furthermore, unit commitment decisions are formulated for each thermal resource cluster using the aggregated representation proposed in~\cite{palmintier_heterogeneous_2014}. {In the context of capacity expansion models, it was observed by~\cite{poncelet_unit_2020} that the relaxation of integrality constraints on clustered unit commitment variables has limited impact on model results for those studies that do not explicitly focus on cycling of thermal generators and their ability to provide flexibility. This is in line with the results in \cite{li_mixed-integer_2022}, showing that the optimality gap obtained by relaxing the integrality constraints is often small in practice. Therefore, this study considers CEMs with continuous unit commitment decision variables, which is also the most common practice in best-in-class open-source capacity expansion models (e.g. PyPSA~\cite{brown_pypsa_2018}, {Switch}~\cite{johnston_switch_2019}, GenX~\cite{mit_energy_initiative_genx_nodate}). We note that when integer unit commitment variables are needed, it is still possible to: (1) relax unit commitment variables to be continuous; (2) apply the Benders decomposition algorithms described in Section~\ref{sec:benders} and~\ref{sec:benders_reg} and run to convergence to obtain a lower bound; (3) fix all planning decisions to the values computed in the previous step and solve the operational sub-problems with integer unit commitment variables obtaining an upper bound (where slacks are used to ensure feasibility).}

The modeling of unit commitment, ramping limits, and storage operation requires the formulation of time-coupling constraints. CEMs often represents them by using {cyclic constraints}, where the first time step of each sub-period is treated as immediately following the last time step of the sub-period for purposes of time-coupling constraints {(for full details, see the documentation of widely used open-source power system capacity expansion models PyPSA~\cite{brown_pypsa_2018}, {Switch}~\cite{johnston_switch_2019}, or GenX~\cite{mit_energy_initiative_genx_nodate}).} As an illustrative example, we discuss this approximation for storage constraints, but the same idea applies to hydropower operation and thermal power plant ramping and unit commitment constraints. {Let $W_{p,s}\subseteq W_p$ be the index set of operational sub-periods of a given scenario $s$ for planning period $p$.} Let $H_w$ be the index set of all hourly time steps included in sub-period $w \in W_{p,s}$, and define ${t^{\text{first}}_w} \in H_w$ as the first time step of the sub-period and ${t^{\text{last}}_w} \in H_w$ as the last time step of the sub-period. Variable $x^{\text{soc}}_{g,t}$ denotes the component of vector $x_w$ corresponding to the state of charge of storage resource $g$ at time $t \in H_w$. Analogously, ${x^{c}_{g,t}}$ and ${x^d_{g,t}}$ correspond to the electricity charged from and discharged into the system, respectively, by storage resource $g$ at time $t \in H_w$. Parameters $\eta^{self}_g$, $\eta^{c}_g$, and $\eta^{d}_g$ are percentages of self-discharge rate, charging efficiency, and discharging efficiency, respectively. All storage resources are modeled by the following equations:
 \begin{equation}
    \label{eq:storage}
    \begin{split}
    &x^{\text{soc}}_{g,t} = (1-\eta^{self}_g)x^{\text{soc}}_{g,t-1} + \eta^{c}_g {x^{c}_{g,t}} - \frac{{x^d_{g,t}}}{\eta^{d}_g}, \\
    &\forall t \in H_w \setminus \{{t^{\text{first}}_w}\}, \quad w \in W_{p,s}\\
    &x^{\text{soc}}_{g,{t^{\text{first}}_w}} = (1-\eta^{self}_g){z^{\text{start}}_{g,w}}  + \eta^{c}_g {x^{c}_{g,t^{\text{first}}_w}} - \frac{{x^d_{g,t^{\text{first}}_w}}}{\eta^{d}_g}, \\
    &\forall w\in W_{p,s},
    \end{split}
\end{equation}
{where $x^{\text{soc}}_{g,t}$ is the storage level (measured as energy) of resource $g$ at time $t$, and $z^{\text{start}}_{g,w}$ is an auxiliary variables representing the storage level at the start of sub-period $w$:
\begin{equation}
    \label{eq:longstor}
    \begin{split}
    &z^{\text{start}}_{g,w+1} =  x^{soc}_{g,t^{\text{last}}_{w}}  \quad \forall w \in W_{p,s}\setminus \{w^{\text{last}}_s\}\\
    &z^{\text{start}}_{w^{\text{first}}_s}= x^{soc}_{g,t^{\text{last}}_{w^{\text{last}}_s}},\\
    \end{split}
\end{equation}
where $w^{\text{first}}_s$ and $w^{\text{last}}_s$ are indices of first and last sub-periods in $W_{p,s}$.
The {cyclic} approximation assumes that storage levels between sufficiently long sub-periods are decoupled. Hence, it discards~\eqref{eq:longstor} and sets:
\begin{equation}
    \label{eq:shortstor}
    z^{\text{start}}_{g,w} = x^{\text{soc}}_{g,{t^{\text{last}}_w}}, \quad \forall w \in W_{p,s}.\\ 
\end{equation}
}
As a result, storage constraints~\eqref{eq:storage} and \eqref{eq:shortstor} corresponding to different sub-periods are independent and are included in~\eqref{eq:opcons}. Such approximation may give rise to errors when sub-periods are too short, as they fail to fully capture electricity demand and variable renewable resources' availability profiles~\cite{mallapragada_impact_2018,jacobson_quantifying_2024}. These errors are exacerbated when considering only few 24-hour sub-periods (days), as done in \cite{Lara2018,li_mixed-integer_2022}. In comparison to existing literature, previous work~\cite{jacobson_computationally_2024} considered $52$ week-long (168 hour) sub-periods for each planning period. While {cyclic constraints} with weekly sub-periods may not cause significant errors when considering short-duration storage technologies~\cite{mallapragada_long-run_2020}, this does not hold for multi-day storage resources (MDS)\cite{mantegna_establishing_2024}. As noted in \cite{levin_energy_2023}, accurately representing multi-day energy storage (MDS) in CEMs is a key requirement to evaluate their role in decarbonized energy systems. {Here, we consider CEMs that accurately model consecutive hourly operation of MDS resources (including hydropower) with constraints~\eqref{eq:longstor}, while all other storage resources, unit commitment, and ramping constraints are modeled using cyclic constraints and week-long sub-periods with hourly timesteps.
To further highlight problem structure, we define a new auxiliary variable $z^{\text{end}}_{g,w}$ to represent storage level at the end of sub-period $w$. Multi-day energy storage constraints~\eqref{eq:longstor} are equivalently rewritten as:
\begin{align}
&z^{\text{end}}_{g,w} = x^{soc}_{g,t^{\text{last}}_{w}} \quad \forall w\in W_{p,s}\label{eq:longstor_decomp} \\
&z^{\text{start}}_{g,w+1} = z^{\text{end}}_{g,w}  \quad \forall w\in W_{p,s}\setminus \{w^{\text{last}}_s\}\label{eq:longstor_planning_1}\\
&z^{\text{start}}_{w^{\text{first}}_s}=z^{\text{end}}_{g,w^{\text{last}}_s}\label{eq:longstor_planning_2}
\end{align}
With this reformulation, constraints~\eqref{eq:longstor_decomp} are grouped in~\eqref{eq:auxcons1}, while~\eqref{eq:longstor_planning_1} and  \eqref{eq:longstor_planning_2} are included in~\eqref{eq:couplingcons}. 
We note that, while a week-long cyclic approximation of unit commitment constraints is adequate to represent operation of fossil burning power plants, more complex fuel cycles of nuclear power plants~\cite{jenkins_benefits_2018} or maintenance planning decisions for thermal plants~\cite{palmintier_impact_2016} could be represented with analogous formulations. In fact, we can included the commitment states linking consecutive sub-periods within vector $z_w$.

Finally, we consider policies like $\text{CO}_2$ emission caps or Renewable Portfolio Requirements (RPS) that link operating decisions across sub-periods. As an example, a $\text{CO}_2$ emission cap enforced for the considered weather realization of period $p \in P$ is:
\begin{equation}
\label{eq:co2_cap}
    \sum_{w \in W_{p,s}}\Bigg(\Big(\sum_{t \in H_w}\sum_{g \in G} \gamma^{\text{CO}_2}_g x^{\text{d}}_{g,t} \Big) \Bigg)\leq \epsilon^{\text{CO}_2},
\end{equation}
where $G$ is the set of all resources, $\gamma^{\text{CO}_2}_g$ represents the tons of $\text{CO}_2$ emissions per-MWh of resource $g \in G$, variable $x^{\text{d}}_{g,t}$ is the electricity that is discharged by resource $g$ at hour $t$, while $\epsilon^{\text{CO}_2}$ is the emission cap. Note that other $\text{CO}_2$ emission cap or minimum energy share requirement formulations can be considered, for example defining the requirement as fraction of total power generation or demand. However, our derivation is analogous for all these alternatives, and we use equation \eqref{eq:co2_cap} for the sake of simplicity. We reformulate policy constraints by decoupling the sub-periods using budgeting variables as proposed in~\cite{jacobson_computationally_2024}. We have shown that constraint~\eqref{eq:co2_cap} is equivalent to:
\begin{align}
    &\sum_{t \in H_w}\sum_{g \in G} \gamma^{\text{CO}_2}_g x^{\text{d}}_{g,t} \leq z^{b}_w, \quad \forall w \in W_p\label{eq:co2_cap_decomp} \\
    &\sum_{w \in W_{p,s}}z^b_w = \epsilon^{\text{CO}_2} \label{eq:co2_cap_planning}\\
    &z^b_w \geq 0
\end{align}
where budgeting variables $z^b_w$ are introduced to decouple the policy constraint on each sub-period. We note that constraints~\eqref{eq:co2_cap_decomp} are represented by~\eqref{eq:auxcons1}, while~\eqref{eq:co2_cap_planning} is included in~\eqref{eq:couplingcons}. 
\section{Benders Decomposition}
\label{sec:benders}
{
{The capacity expansion models described in Section~\ref{sec:multi_cems} can be obtained with widely used open-source power system capacity expansion models PyPSA~\cite{brown_pypsa_2018}, {Switch}~\cite{johnston_switch_2019}, or GenX\cite{mit_energy_initiative_genx_nodate}. Furthermore, the majority of CEMs considered in previous academic literature can be written as Problem~\eqref{eq:probform}, including \cite{munoz_new_2016,Lara2018,li_mixed-integer_2022,zhang_integrated_2023,soares_integrated_2022,jacobson_computationally_2024,goke_stabilized_2024}. }
} However, due to the computational challenges encountered when solving models with full-year resolution, existing capacity expansion models are often simplified by considering only representative sub-periods, ignoring MDS resources, and/or relaxing integrality constraints for investment and retirement decisions.
Previous work~\cite{jacobson_computationally_2024} focused only on a Benders decomposition approach for policy constraints and did not consider MDS resources, which introduce additional time-coupling constraints~\eqref{eq:longstor}. Here, we move beyond existing literature and allow sub-period decomposition of the operational year for CEMs that include MDS constraints with hourly resolution and time-coupling policy constraints.}
Denote by $y$ and ${z}$ the vectors consisting of $y_p$ and ${z_w}$ for all $w \in W_p$ and $p \in P$, respectively.
At iteration $k \geq 0$, given a choice of planning decisions $y^{(k)}$ and ${z^{(k)}}$, we solve the following operational sub-problem for each $w\in W_p$ and $p \in P$:
\begin{subequations}
    \label{eq:subprob}
    \begin{alignat}{3}
        g_w^{(k)} =&\text{min}&\; \; &c_w^T x_w \label{eq:obj_subprob}\\
        &\text{s.t.}& & A_wx_w + B_wy_p \leq b_w \label{eq:opcons_subprob}\\
        &&& Q_w x_w \leq {z_w}\\
        &&& y_p =y_p^{(k)} \quad {(\pi_w^{(k)})} \\
        &&& {z_w = z_w^{(k)}} \quad {(\lambda_w^{(k)})} \\
        &&& x_w \geq 0,
    \end{alignat}
\end{subequations}
where ${\pi_w^{(k)}}$ and $\lambda^{(k)}_w$ are Lagrangian multipliers associated with the corresponding constraints. We define the current best upper bound as:
\begin{equation}
    \label{eq:upper_bound_comp}
    U^{(k)} = \min_{j=0,\ldots,k}\bigg(\sum_{p \in P}\Big(f_p^Ty_p^{(j)} + \sum_{w \in W_p}{\beta_{\sigma(w)}}g_w^{(j)}\Big)\bigg),
\end{equation}
and denote the planning decisions responsible for the best upper bound as $y^*$ and $z^*$.
Next, we compute new planning decisions $y^{(k+1)}$ and ${z^{(k+1)}}$ by solving the {planning} problem:
\begin{subequations}
    \label{eq:planning}
    \begin{alignat}{3}
        L^{(k)}=&\text{min}&\; &\sum_{p \in P}\Big(f_p^Ty_p+\sum_{w \in W_p}{\beta_{\sigma(w)}}\theta_w\Big) \label{eq:obj_planning}\\
        &\text{s.t.}& & \theta_w \geq g_w^{(j)} + (y_p-y_p^{(j)})^T{\pi_w^{(j)}} + ({z_w-z_w^{(j)}})^T\lambda_w^{(j)},\label{eq:bdcuts_planning}\\ &&& \quad \forall j=0,\ldots,k, w \in W_p, \; p \in P \nonumber \\
        &&& C_w z_w + D_wy_p \leq d_w, \quad \forall w \in W_p, \; p \in P \\
        &&& {\sum_{w \in W_{p,s}}E_wz_w = e_{p,s}, \quad \forall s \in S,\;\forall p \in P} \label{eq:couplingcons_planning}\\
        &&& \sum_{p \in P}R_p y_p \leq r \label{eq:multiperiodcons_planning} \\
        &&& z_w \geq 0, \quad \forall w \in W_p, \; p \in P\\
        &&& y_p \geq 0,\quad \forall p \in P\\
        &&& y_{p,i} \in \mathbb{Z}, \quad \forall i \in I, \; \forall p \in P.\label{eq:int_invcons}
    \end{alignat}
\end{subequations}
and proceed to the next iteration. The resulting Benders decomposition algorithm is presented in Algorithm~\ref{alg:benders}.
\begin{algorithm}
    \caption{Benders decomposition.}
    \label{alg:benders}
    \begin{algorithmic}
        \State \textbf{Input:} Set $K_{\max}$, $\epsilon_{\text{tol}}$, and $\alpha \in (0,1)$.
        \State \textbf{Output:} $y^*$ and ${z^*}$.
        \State { \textbf{Initialization:} Solve Problem~\eqref{eq:planning} without any Benders cuts and obtain  $y^0$ and ${z^0}$. }
        \For{$k=0,\ldots,K_{\max}$}
        \For{$p \in P$}
        \For{$w \in W_p$}
        \State Solve operational sub-problem~\eqref{eq:subprob}.
        \EndFor
        \EndFor
        \State{Compute upper bound $U^{(k)}$ and update $y^*$ and $z^*$.}
        \State {Update cuts in Problem~\eqref{eq:planning} and solve it to get $L^{(k)}$}.
        \If{$(U^{(k)}-L^{(k)})/L^{(k)} \leq \epsilon_{\text{tol}}$} 
        \State Stop and return {$y^*$ and $z^*$}.
        \Else{$\;$ Set $y^{(k+1)}$ and ${z^{(k+1)}}$ equal to the solution of~\eqref{eq:planning}.}
        \EndIf
        \EndFor 
    \end{algorithmic}
\end{algorithm}
Compared to previous literature, each iteration solves multiple smaller sub-problems~\eqref{eq:subprob} formulated over sub-periods of the operational year, adding more Benders cuts per iteration. As shown in~\cite{jacobson_computationally_2024} this greatly improves convergence compared to solving a full size sub-problem for each planning period, which would add one cut per period at each iteration, as done in~\cite{Lara2018,li_mixed-integer_2022,zhang_integrated_2023,goke_stabilized_2024} - see also the results in Section~\ref{sec:numexp}. Moreover, the sub-problems can be dealt with in parallel, significantly reducing the wall clock time per iteration compared to solving a full size sub-problem. Algorithm~\ref{alg:benders} extends previous work~\cite{jacobson_computationally_2024} by considering multi-period investment and retirement decisions, rather than a single-period CEM, as well as MDS resources. Note that we assume that sub-problems~\eqref{eq:subprob} are always feasible, for every choice of $y_p^{(k)}$ and ${z_w^{(k)}}$. {To this purpose, we include slack variables and penalty terms both for policy constraints~\eqref{eq:co2_cap_decomp} and multi-day storage constraints~\eqref{eq:longstor_decomp}. In general, the choice of these penalty terms will depend on the focus of the specific study. For example, if interested in investigating the role of multi-day energy storage resources, modelers could choose a large penalty to enforce strict satisfaction of the constraints. As an alternative to introducing slacks, Algorithm\eqref{alg:benders} could include the computation of feasibility cuts when the operational sub-problem is infeasible - as example, see~\cite{kronqvist_using_2020}. This option comes with its own challenges associated with the definition of tight feasibility cuts, but it does not fundamentally modify the nature of the algorithms described here}.

\section{Regularized Benders Decomposition}
\label{sec:benders_reg}
Algorithm~\ref{alg:benders} belongs to the class of cutting plane methods, and it is expected to suffer from instability\cite{ruszczynski_regularized_1986,lemarechal_new_1995}.
As shown in Section~\ref{sec:numexp}, this behavior affects Algorithm~\ref{alg:benders} when implemented as done in \cite{jacobson_computationally_2024}, without any regularization technique. In fact, Algorithm~\ref{alg:benders} may oscillate between different extreme investment and retirement decisions before being able to make substantial progress towards the solution. This is because the approximated system cost in~\eqref{eq:obj_planning} may significantly underestimate the operational cost of extreme decisions when the number of cuts~\eqref{eq:bdcuts_planning} is small compared to the number of variables.
{Here, we move beyond previous literature~\cite{jacobson_computationally_2024} and show that level-set methods~\cite{lemarechal_new_1995,zhang_integrated_2023,kronqvist_using_2020} can regularize Algorithm~\ref{alg:benders} enabling high-performance capacity expansion models with hourly resolution.}
Let $\Phi(\cdot)$ be a convex function and $\alpha \in (0,1)$. At iteration $k\geq 1$ of the Benders decomposition algorithm, we consider the following regularization problem:
\begin{subequations}
    \label{eq:reg_problem}
    \begin{alignat}{3}
        &\text{min}&\; & \Phi(y,{z}) \\
        &\text{s.t.}& & \theta_w \geq g_w^{(j)} + (y_p-y_p^{(j)})^T{\pi_w^{(j)}} + {(z_w-z_w^{(j)}})^T\lambda_w^{(j)},\\ 
        &&& \quad \forall j=0,\ldots,k, w \in W_p, \; p \in P \nonumber \\
        &&& C_w z_w + D_wy_p \leq d_w, \quad \forall w \in W_p, \; p \in P \\
        &&&  {\sum_{w \in W_{p,s}}E_wz_w = e_{p,s}, \quad \forall s \in S,\;\forall p \in P}\\
        &&& \sum_{p \in P}R_p y_p \leq r  \\
        &&& \sum_{p \in P}\Big(f_p^Ty_p+\sum_{w \in W_p}\theta_w\Big) \leq L_{\alpha}^{(k)}\label{eq:level-set-con}\\
        &&& y_p \geq 0,\quad \forall p \in P\\
        &&& y_{p,i} \in \mathbb{Z}, \quad \forall i \in I, \; \forall p \in P.
    \end{alignat}
\end{subequations}
where $L_{\alpha}^{(k)} = L^{(k)}+\alpha(U^{(k)}-L^{(k)})$. In contrast to Algorithm~\ref{alg:benders}, we do not define the next iterate $y^{(k+1)}$ and ${z^{(k+1)}}$, as the solution of Problem~\eqref{eq:planning}, but rather as the solution of Problem~\eqref{eq:reg_problem}, which selects the best feasible solution according to the criterion $\Phi(\cdot)$, and subject to the level-set constraint~\eqref{eq:level-set-con} on the approximated system cost. { Note that the solution of the planning Problem~\eqref{eq:planning} belongs to the the feasible set of Problem~\eqref{eq:reg_problem}. In particular, it satisfies the level-set constraint~\eqref{eq:level-set-con}. Hence, as long as the planning Problem~\eqref{eq:planning} is feasible, so is the regularization Problem~\eqref{eq:reg_problem}.}
\begin{algorithm}
    \caption{Regularized Benders decomposition.}
    \label{alg:reg_benders}
    \begin{algorithmic}
        \State \textbf{Input:}Set $K_{\max}$, $\epsilon_{\text{tol}}$, and $\alpha \in (0,1)$.
        \State \textbf{Output:} $y^*$ and ${z^*}$.
        \State { \textbf{Initialization:} Solve Problem~\eqref{eq:planning} without any Benders cuts and obtain  $y^0$ and ${z^0}$. }
        \For{$k=0,\ldots,K_{\max}$}
        \For{$p \in P$}
        \For{$w \in W_p$}
        \State Solve operational sub-problem~\eqref{eq:subprob}.
        \EndFor
        \EndFor
        \State{Compute upper bound $U^{(k)}$ and update $y^*$ and $z^*$.}
        \State {Update cuts in Problem~\eqref{eq:planning} and solve it to get $L^{(k)}$}.
        \If{$(U^{(k)}-L^{(k)})/L^{(k)} \leq \epsilon_{\text{tol}}$} 
        \State Stop and return {$y^*$ and $z^*$}.
        \Else{$\;$ Solve~\eqref{eq:reg_problem} to get $y^{(k+1)}$ and ${z^{(k+1)}}$}.
        \EndIf
        \EndFor 
    \end{algorithmic}
\end{algorithm}

Previous literature~\cite{lemarechal_new_1995,kronqvist_using_2020} have suggested various choices for $\Phi(y,z)$, with the most common being $\Phi^{\ell_2}(y,{z}) = \| y-{y^{*}} \|^2_2 + \|{z}-{z^{*}}\|^2_2$,
which turns Problem~\eqref{eq:reg_problem} into a quadratic optimization problem, while the planning problem{~\eqref{eq:planning}} is linear. {Hence, this strategy is beneficial if it can compensate for the increased computational effort by reducing the number of total iterations of Algorithm~\ref{alg:benders}.}
As reported in Table~\ref{tab:results}, this approach did not perform well when applied to the largest case study considered in Section~\ref{sec:numexp}, both in terms of runtime and number of iterations. An alternative method adopted by \cite{zhang_integrated_2023} aims to select solutions that are $\alpha(U^{(k)}-L^{(k)})$-suboptimal and belong to the interior of the feasible set of Problem~\eqref{eq:planning}. This strategy is frequently used in decomposition algorithms when solving continuous problems~\cite{gondzio_new_2013}. Because integer variables take values at the bounds and do not belong to the interior of the feasible set, this regularization method can not be directly applied to MILPs. For continuous optimization problems, an interior solution can be computed by setting $\Phi(y,{z}) =\Phi^{\text{int}}(y,{z}) = 0$, and solving the resulting feasibility problem with a barrier method. {In section~\ref{sec:numexp}, we use the solver Gurobi with the barrier algorithm~(\texttt{Method}=2) and without crossover~(\texttt{Crossover}=0). 

To benefit from interior-point regularization when MILPs are considered, we design a two-stage method. First implement Algorithm~\ref{alg:reg_benders} with $\Phi=\Phi^{\text{int}}$ to solve the continuous relaxation of Problem~\eqref{eq:probform}. Note that, even though they have been computed using continuous planning decisions, the cuts generated by Algorithm~\ref{alg:reg_benders} do not exclude any integer optimal solution of Problem~\eqref{eq:probform}. Therefore, they can be used to warm-start the planning problem when solving the original MILP. To recover an integer solution, we then apply Algorithm~\ref{alg:reg_benders} keeping all pre-computed cuts in Problem~\eqref{eq:planning} (which now includes integer constraints~\eqref{eq:int_invcons}) and using a modified regularization sub-problem. In this case, we fix discrete investment and retirement decisions to the values computed solving the planning problem and select the remaining continuous planning decisions in the interior of the feasible set.
 The two-stage method is summarized in Algorithm~\ref{alg:reg_benders_int}.
 } 

\begin{algorithm}
    \caption{Benders-$\Phi^{\text{int}}$ decomposition.}
    \label{alg:reg_benders_int}
    \begin{algorithmic}
        \State \textbf{Stage 1.} Apply Algorithm~\ref{alg:reg_benders} with $\Phi=\Phi^{\text{int}}$ ignoring integrality constraints~\eqref{eq:int_invcons}. 
        \State       
        \State \textbf{Stage. 2} Apply Algorithm~\ref{alg:reg_benders} including all integrality constraints and initializing Problems~\eqref{eq:planning} and~\eqref{eq:reg_problem} with the cuts computed by Stage 1. At the end of each iteration, solve Problem~\eqref{eq:reg_problem} with additional constraints: 
        \begin{equation}
        y_{p,i} = \hat{y}_{p,i}^{k+1}, \; \forall i \in I, \; p \in P,
        \end{equation}
        where $\hat{y}_p^{k+1}$ is the vector of integer planning decisions found solving Problem~\eqref{eq:planning}.
    \end{algorithmic}
\end{algorithm}

\section{Results}
\label{sec:numexp}
In this section, we numerically evaluate the proposed algorithms and compared them to previously published alternatives. {The Benders decomposition framework described in Sections~\ref{sec:benders} and~\ref{sec:benders_reg} applies to all CEMs that results in optimization problems like Problem~\eqref{eq:probform}. For example, these can be obtained using open-source modeling packages PyPSA~\cite{brown_pypsa_2018} or {Switch}~\cite{johnston_switch_2019}. In this study, we utilize the open-source package GenX~(\href{https://github.com/GenXProject/GenX/releases/tag/v0.3.5}{v0.3.5})\footnote{All data and codes used in these numerical experiments are available at \url{https://doi.org/10.5281/zenodo.12724093}.}, with full model details available in the documentation~\cite{mit_energy_initiative_genx_nodate}. We implemented our Benders decomposition algorithms as new solving routines within GenX.} All LPs, MILPs, and QPs are solved using the solver Gurobi~(v.10.0.1), and the GenX model is implemented in Julia (v1.9.1)\cite{bezanson_julia_2017} and JuMP (v1.17)\cite{lubin_jump_2023}. {All LPs and QPs are solved with the barrier method, with option \texttt{BarConvTol} set to $10^{-3}$. For MILPs, we also set option \texttt{MIPGap} to $10^{-3}$.} Analogously, we set the convergence tolerance of Algorithms~\ref{alg:benders},~\ref{alg:reg_benders}, and~\ref{alg:reg_benders_int} to $10^{-3}$. When solving~\eqref{eq:planning} and~\eqref{eq:reg_problem} we deactivate crossover, which is switched on to compute optimal basic primal-dual solutions of the operational sub-problems~\eqref{eq:subprob}. {Non-basic dual feasible solutions of the operational sub-problems computed by the barrier method without performing crossover would still provide valid Benders cuts. However, these cuts may not separate the optimal solution from the current primal solution if the convergence tolerance is not small enough~\cite{zakeri_inexact_2000}. The use of cuts obtained from non-basic solutions may increase the number of iterations and it requires to progressively reduce the convergence tolerance for the sub-problems in order to guarantee convergence~\cite{zakeri_inexact_2000}. Hence, we use crossover to obtain better quality cuts. The simulations are run on Princeton University’s Della computer cluster, constraining all experiments to run on 2.8 GHz Intel Cascade Lake nodes\footnote{\url{https://researchcomputing.princeton.edu/systems/della\#hardware}}. The distribution and solution of operational sub-problems on different cores is done through the Julia package \texttt{Distributed.jl}. To avoid issues with shared memory access, we set the Gurobi parameter \texttt{Threads}=1 for all operational sub-problems.}
\subsection{Multi-period capacity expansion with discrete decisions}
\begin{figure}
    \centering
    \subfloat[3 zones\label{fig:conus_3z}]{%
      \includegraphics[width=0.20\textwidth]{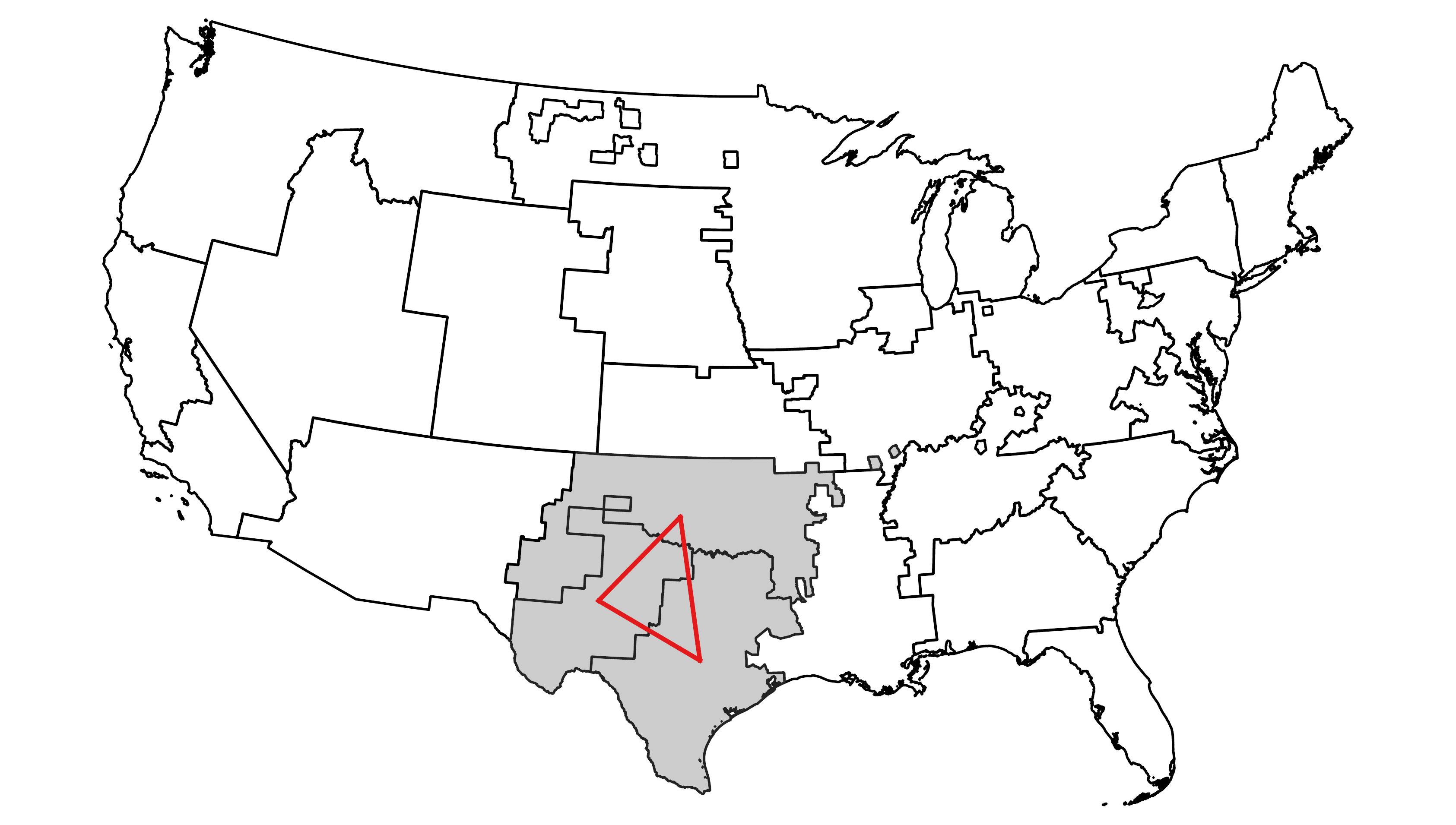}
    }
    \subfloat[6 zones\label{fig:conus_6z}]{%
      \includegraphics[width=0.20\textwidth]{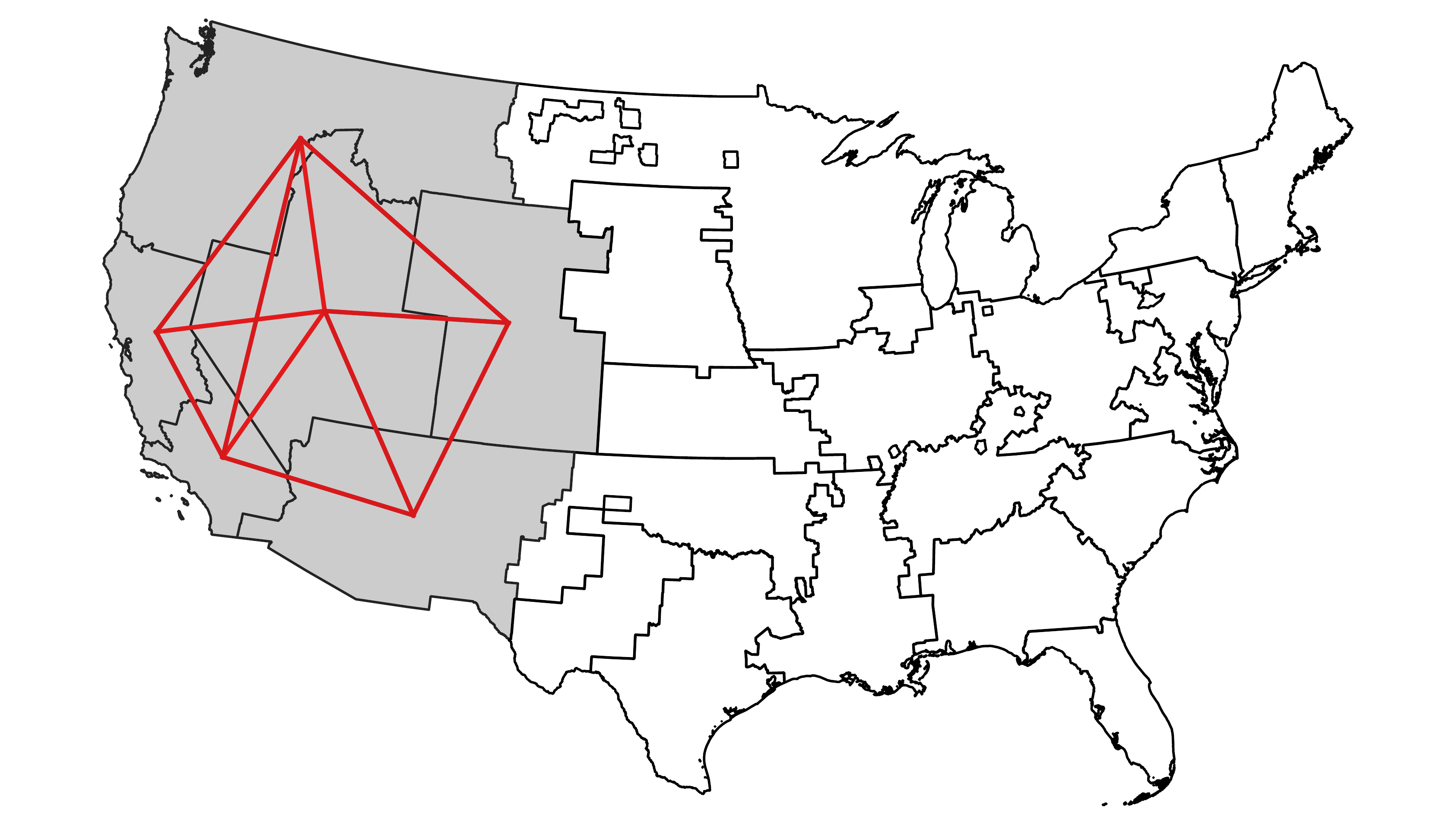}
    }
    \hfill
    \subfloat[15 zones\label{fig:conus_15z}]{%
      \includegraphics[width=0.20\textwidth]{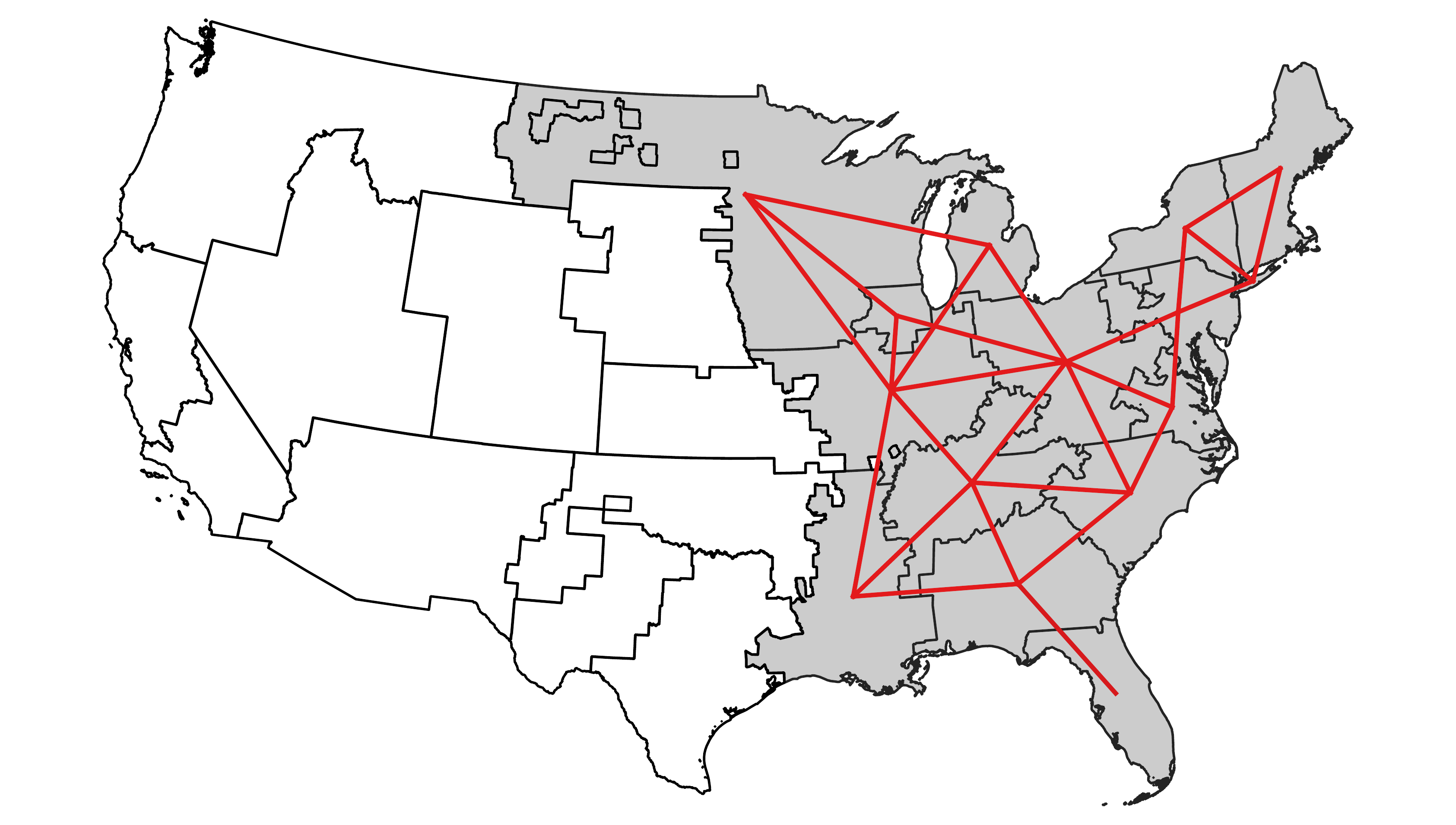}
    }
    \subfloat[26 zones\label{fig:conus_26z}]{%
      \includegraphics[width=0.20\textwidth]{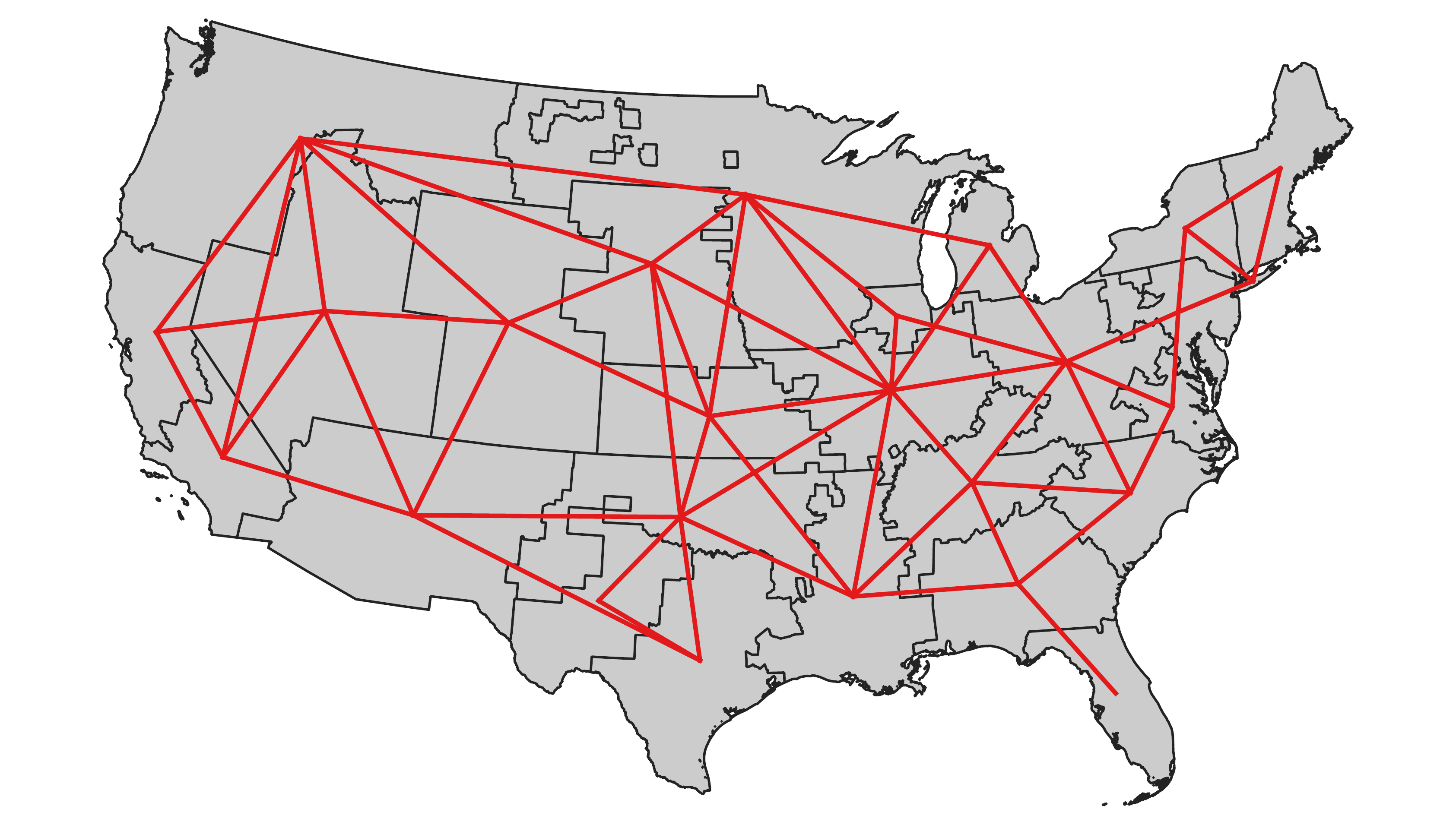}
    }
    \caption{Test systems from the Continental United States (CONUS).}
    \label{fig:system_map}
\end{figure}

We consider four electricity system models from the Continental United States with input parameters generated by PowerGenome\footnote{\url{https://github.com/PowerGenome/PowerGenome}}- see Figure~\ref{fig:system_map}. {The largest system with 26 zones has 1004 storage and generation resources, including 408 thermal generator clusters, 431 variable renewable energy clusters, $16$ hydropower reservoirs, and $24$ hydroelectric pumped storage. Both hydropower reservoirs and hydroelectric pumped storage resources are modeled as MDS.} In addition, the 26 zone electricity network includes 49 transmission links representing existing power transmission capacity between zones, and 8 candidate links corresponding to possible new connections. Both existing and candidate transmission links are eligible for expansion. {Discrete investment decisions explicitly model the number and type of lines that are built on each connection. For each of the inter-zonal connections in Figure~\ref{fig:system_map} we consider three line voltages (230kV, 345kV, 500kV) with either single or double circuits, and a 500kV HVDC line, resulting in 7 different line classes. Discrete transmission expansion decisions correspond to the number of new lines from each class. For the largest 26-zone system, this results in $57\times 7=399$ transmission expansion integer variables for each planning period. Analogously, discrete generation investment or retirement decisions model the number of units in generation and storage resource clusters that are built or retired.}

We consider transmission and generation capacity expansion over 3 planning periods (2024-2030, 2031-2040, 2041-2050), with each planning stage represented by a single operational year (8736 hours), with demand corresponding to the final year of the planning stage (to ensure resource adequacy) and investment costs corresponding to the average costs over the planning stage (representing the average cost of a series of capacity additions across the planning stage). An annual CO$_2$ emission constraint is enforced on the operational decisions corresponding to each stage individually. The CO$_2$ cap in each year is set to { represent a pathway to net-zero by 2050}, with values of 186 Mt for 2030, 86.66 Mt for 2040 and 0 Mt for 2050. 
{To ensure feasibility of the sub-problems we add slacks to the CO$_2$ emission constraints, and set a penalty value for violating these constraints equal to $150\$/\text{ton}$. Depending on the focus of the study, larger CO$_2$ price values can be used to limit violations of the CO$_2$ cap.
Moreover, the penalty for violating multi-day storage constraint at the beginning and end of each subperiod is set equal to twice the cost of non-served energy demand. In this way, sub-problems~\eqref{eq:subprob} are always feasible and we ensure that MDS constraints are met in the final optimal solution before any non-served demand.} System operation is modeled with hourly resolution, and subdivided into 52 sub-periods (i.e. weeks) for each planning period. Therefore, each sub-period includes 168 time steps and the total number of operational time steps is $26,208$. The resulting capacity expansion model for the largest system with 26 zones is written as Problem~\eqref{eq:probform} {with over $111$ million variables, $197$ million constraints, and $4,650$ integer variables}. Note that the considered model has significantly higher spatial and technological resolution compared to those considered in previous literature~\cite{Lara2018,li_mixed-integer_2022,zhang_integrated_2023,goke_stabilized_2024}, resulting in a capacity expansion model that has from 2 to 20 times more variables.

First, we compare the performance of several solution algorithms to solve the continuous relaxation of the multi-period planning Problem~\eqref{eq:probform}, where the integrality constraints in~\eqref{eq:int_invcons} are ignored. Algorithm~\ref{alg:reg_benders}, referred to as \texttt{Benders-}$\Phi^{\alpha}$, is implemented with $\Phi \in \{\Phi^{\text{int}},\Phi^{\ell_2}\}$ and level-set parameter $\alpha \in \{0.2,0.5,0.8\}$. To assess the effect of regularization on the convergence properties of the Benders algorithm, we also evaluate Algorithm~\ref{alg:benders}, referred to as \texttt{Benders-NoReg}. In addition, we performed benchmarking experiments for algorithms proposed in previous literature. These include the Dual Dynamic Programming (\texttt{DDP}) decomposition algorithm proposed in \cite{Lara2018} and previously implemented in GenX (as of v0.3.0), which does not allow any parallel implementation. {In fact, DDP decomposes our multi-period CEM into three single-period CEMs, which have to be solved in series at each iteration.} 
Furthermore, the trust-region scheme with $\ell_{\infty}$-norm from \cite{goke_stabilized_2024}, denoted by \texttt{Benders-TR}, is applied as regularization step in Algorithm~\ref{alg:reg_benders}. Because it was observed in \cite{goke_stabilized_2024} that the algorithm performance did not significantly depend on its parameters, we set most parameters as in \cite{goke_stabilized_2024} and only vary the trust-region parameter in $\{0.1,0.2,0.3\}$. 

{The results obtained by regularized Benders algorithms for different parameter values are reported in full in the public repository$^1$. Table~\ref{tab:results} reports only the performance of the parameter choice that resulted in either the smallest number of iterations to converge or the lowest remaining gap for each regularization method. Note that \texttt{Benders-}$\Phi^{\alpha}$, \texttt{Benders-TR}, and \texttt{Benders-noReg} solve $156$, parallelized, week-long, operational sub-problems at each iteration. In comparison, the \texttt{DDP} algorithm from \cite{Lara2018} is required to solve six year-long, single-period, planning problems in each iteration to compute upper and lower bounds, as it must proceed through forward and backward passes in series. Hence, \texttt{DDP} does not allow parallelization, limiting its scalability performance with the number of available CPU cores. We expect the individual iterations of \texttt{DDP} to require significantly more computational effort than those of the Benders-based method. }

All \texttt{Benders-}$\Phi^{\alpha}$, \texttt{Benders-TR}, and \texttt{Benders-noReg} experiments reported in this section {use 156 cores distributed over 6 computing nodes, allowing complete parallelization of the 156 operational sub-problems.} 
In comparison, \texttt{DDP}~\cite{Lara2018} is implemented on a single node with {96 cores (no parallel processes)}. {For each algorithm, we report computational time needed to actually solve the model, without considering input/output routines and model generation, and we set a time limit of 12 hours. \texttt{DDP} experienced numerical issues when solving the 6-,15- and 26-zone cases, resulting in early termination due to negative optimality gaps. These issues were resolved by reducing Gurobi's parameter \texttt{BarConvTol} to $10^{-8}$ for these cases. Because of the reduced tolerance, we set the time limit for DDP to be 24 hours.}

\begin{table}[!h]
    \caption{Results for the continuous relaxation of Problem~\eqref{eq:probform}. We denote runtime (in minutes) with $t$ and define $k^*$ and $t^*$ as the first iteration and time when the gap goes below 1\%.}
    \label{tab:results}
    \centering
\begin{tabular}{|l|r|r|r|r|r|}
    \hline
    \text{Method} & \text{Gap (\%)} & \text{Iter} & t (60s) & $k^*$ & $\text{t}^*$ (60s)\\\hline
\rowcolor{lightgray}    \multicolumn{6}{|c|}{{3 zones (CEM with over 9M variables and 18M constraints)}}\\\hline
$\texttt{DDP}$ & 0.10 & 4 & 145.77 &   2 & 77.30 \\\hline
$\texttt{Benders-noReg}$ & 0.09 & 30 & 3.38 &   20 & 2.39 \\\hline
$\texttt{Benders-}\Phi^{{\ell_2},0.2}$ & 0.08 & 28 & 3.58 &   17 & 1.78\\\hline
$\texttt{Benders-TR}^{0.3}$ & 0.07 & 64 & 8.40 &   54 & 6.92 \\\hline
\textbf{$\texttt{Benders-}\Phi^{{int},0.5}$} & \textbf{0.07} & \textbf{19} & \textbf{2.10} &  \textbf{ 13} & \textbf{1.49} \\\hline

\rowcolor{lightgray}    \multicolumn{6}{|c|}{{6 zones (CEM with over 24M variables and 43M constraints)}}\\\hline

$\texttt{DDP}$ & 0.09 & 2 & 295.15 &   2 & 295.15 \\\hline
$\texttt{Benders-noReg}$ & 0.09 & 73 & 47.26 &   42 & 26.54 \\\hline
$\texttt{Benders-}\Phi^{{\ell_2},0.2}$ & 0.09 & 68 & 83.20 &   38 & 20.77 \\\hline
$\texttt{Benders-TR}^{0.3}$ & 0.09 & 89 & 58.76 &   74 & 45.95 \\\hline
\textbf{$\texttt{Benders-}\Phi^{{int},0.5}$} & \textbf{0.09} & \textbf{28} & \textbf{17.04} &  \textbf{ 17} & \textbf{9.17} \\\hline

\rowcolor{lightgray}    \multicolumn{6}{|c|}{{15 zones (CEM with over 59M variables and 108M constraints)}}\\\hline

$\texttt{DDP}$ & 15.75 & 1 & 1400.00 & Inf & Inf \\\hline
$\texttt{Benders-noReg}$ & 0.10 & 128 & 504.96 &   65 & 322.03 \\\hline
$\texttt{Benders-}\Phi^{{\ell_2},0.2}$ & 0.26 & 103 & 720.00 &   49 & 278.65 \\\hline
$\texttt{Benders-TR}^{0.3}$ & 0.08 & 89 & 284.54 &   79 & 227.05 \\\hline
\textbf{$\texttt{Benders-}\Phi^{{int},0.5}$} & \textbf{0.09} & \textbf{42} & \textbf{81.88} &  \textbf{ 21} & \textbf{40.28} \\\hline

\rowcolor{lightgray}\multicolumn{6}{|c|}{{26 zones (CEM with over 111M variables and 197M constraints)}}\\\hline

$\texttt{DDP}$ & Inf & 0 & 1440.00 & Inf & Inf \\\hline
$\texttt{Benders-noReg}$ & 234.66 & 11 & 720.00 &  Inf & Inf \\\hline
$\texttt{Benders-}\Phi^{{\ell_2},0.5}$ & 6.42 & 31 & 720.00 &  Inf & Inf \\\hline
$\texttt{Benders-TR}^{0.2}$ & 159.69 & 43 & 720.00 &  Inf & Inf \\\hline
\textbf{$\texttt{Benders-}\Phi^{{int},0.5}$} & \textbf{0.08} & \textbf{51} & \textbf{371.06} &  \textbf{ 27} & \textbf{210.78} \\\hline
    \end{tabular}
    \end{table}

{As shown in Table~\ref{tab:results}, the regularized Benders algorithm \texttt{Benders-}$\Phi^{\text{int},0.5}$ resulted in the best performance across all four testing systems. Not only does \texttt{Benders-}$\Phi^{\text{int},0.5}$ exhibit the lowest runtime among the tested Benders-based methods, but it also requires the smallest number of iterations to converge to a MIP gap of 0.1\%. \texttt{Benders-}$\Phi^{\text{int},0.5}$ is also the quickest algorithm to reach a $1\%$ optimality gap in all tests. Furthermore, the interior-point regularized decomposition method was able to solve the largest case study (26 zones) in a little over 6 hours, while all other tested approaches failed to converge within the time limit. In particular, they all stopped with large optimality gaps, the smallest being $6.42\%$ achieved by \texttt{Benders-}$\Phi^{\ell_2,0.2}$.  As expected, \texttt{DDP} was the slowest of the implemented algorithms. In the 6-, 15- and 26-zone cases, this could be explained by the need to tighten Gurobi's barrier convergence tolerance to avoid numerical instabilities experienced by the \texttt{DDP} method, which increases the number of Gurobi iterations. However, we observe that even for the 3-zone case study, where it was implemented with the same barrier convergence tolerance as the others, \texttt{DDP} had the longest runtime by an order of magnitude. This suggests that \texttt{DDP}'s biggest limitation is the inability to parallelize the solution of the sub-problems. When considering the 15-zone system \texttt{DDP} reached the time limit of 24 hours while performing its second iteration, terminating with optimality gaps larger than $10\%$. In the case of the 26-zone system, \texttt{DDP} did not complete the first iteration before reaching the time limit.} Moreover, when $\texttt{Benders-}\Phi^{\text{int},0.5}$ was applied to the 26-zone case, the computer cluster reported a memory usage {(for both model generation and solution)} roughly equal to {138 GB on each computing node}, compared to the $\sim$447 GB on a single node needed by \texttt{DDP}, which can not exploit distributed computing resources. {Hence, we note that while Benders-based methods use more total memory to setup models on every distributed process, this approach requires allocating less memory on each computing node, which may make it easier to requisition compute resources.}
{We note that runtime of $\texttt{Benders-}\Phi^{\text{int},0.5}$ still grows quadratically with the size of the considered system (e.g., the number of zones and resource options). Hence, future work should investigate the combination of our temporal decomposition across sub-periods with network-based decomposition techniques to increase scalability with respect to the spatial/network dimension.}

\begin{figure}[h!]
    \centering
    \includegraphics*[width=0.40\textwidth]{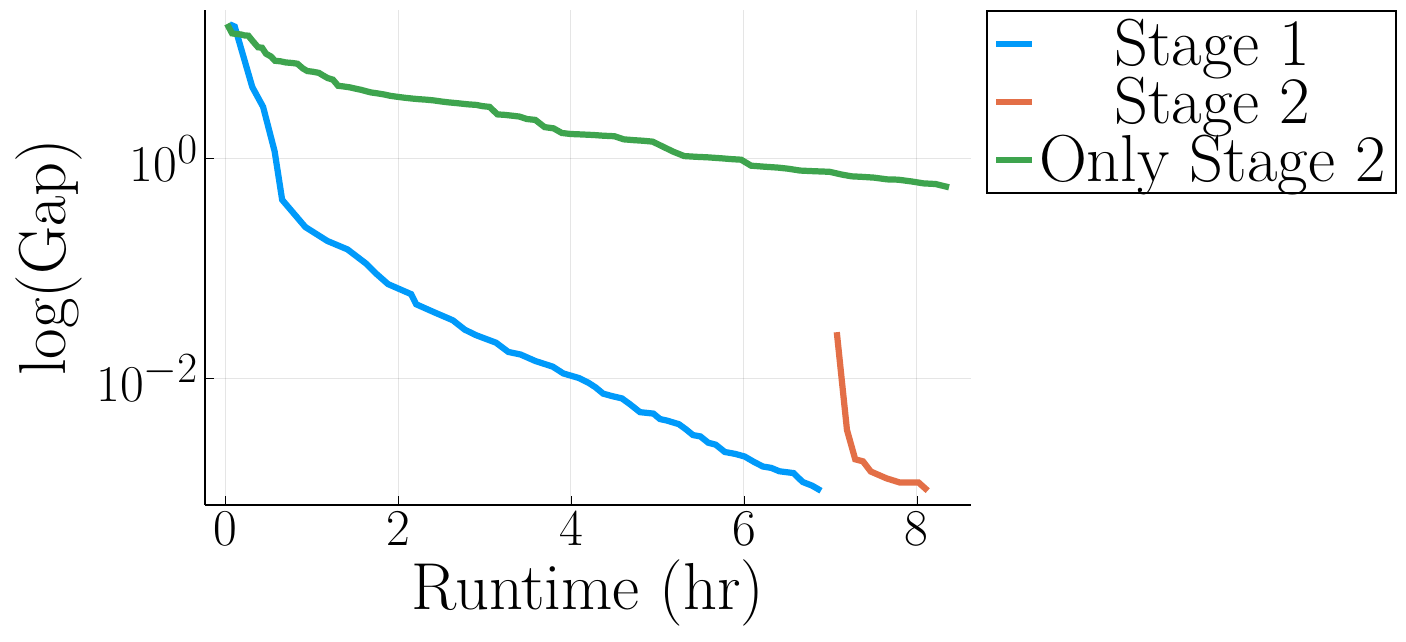}
    \caption{Convergence of Algorithm~\ref{alg:reg_benders_int}, where Stage 1 solves the continuous relaxation, while Stage 2 enforces the integrality constraints on investment and retirement decisions. Note that the y-axis uses a $\log_{10}$ scale.}
    \label{fig:gap_plot}
\end{figure}

To demonstrate the ability to solve CEMs with discrete investment and retirement decisions, we consider the largest test system with 26 zones and enforce the integrality constraint in~\eqref{eq:int_invcons} on investment and retirement decisions. We use Algorithm~\ref{alg:reg_benders_int} with $\alpha=0.5$, which is the regularization scheme that performed best in the previous Section. Recall that Algorithm~\ref{alg:reg_benders_int} performs two stages: the first stage ignores integrality constraints and applies the regularized Benders Algorithm~\ref{alg:reg_benders}. The second stage enforces integrality constraints and implements Algorithm~\ref{alg:reg_benders} where interior-point regularization is applied only to continuous planning variables and Problems~\eqref{eq:planning} and~\eqref{eq:reg_problem} are initialized with the cuts computed at the previous stage. As shown in Figure~\ref{fig:gap_plot} the initialized planning problem at the start of Stage 2 results in a small residual optimality gap ($\sim 1\%$) when integer constraints are taken into consideration. After just four iterations, the optimality gap decreases to below $0.1\%$ reaching convergence. Thus, Algorithm~\ref{alg:reg_benders_int} was able to solve the model in just over 8 hours, {an unprecedented computational performance for a multi-period capacity expansion model applied to an electricity system of this size, modeled with full hourly resolution, and with this number of discrete generation and transmission planning decisions.}
{Finally, we note that the relative gap between best lower bounds of the MILP problem and its continuous relaxation is $0.13\%$. This confirms that the MIP gap for the considered model is small, but still larger than then convergence tolerance of $0.1\%$.} {Figure~\ref{fig:gap_plot} also reports the performance of Algorithm~\ref{alg:reg_benders_int} when we ignore Stage 1 and apply the second stage directly to the MILP problem, considering the integrality constraints from the start. These results show that applying Algorithm~\ref{alg:reg_benders} to the continuous relaxation as a warm-starting procedure is critical to enable the convergence of the decomposition method.}

{Given optimal discrete planning decisions, we can fix their values in Problem~\eqref{eq:probform} and solve it to optimize system’s operation. This allows the evaluation of operational quantities like CO$_2$ emissions. The optimized emissions for the considered test system are equal to 195Mt for 2030, 271Mt for 2040, and 287Mt for 2050. These values are larger than the corresponding CO$_2$ caps, because those constraints are relaxed with slack variables. In this case, we have set a CO$_2$ price equal to 150 \$/ton, which can be interpreted as the price at which a CO$_2$ cap is waived or emission reduction must happen outside of the power sector. Larger penalty values can be selected when strict enforcement of CO$_2$ caps is required. For example, using a penalty equal to twice the cost of non-served demand has enforced satisfaction of the multi-day storage constraints at start and end of sub-periods, with the largest violation being smaller than $10^{-6}$ MWh. 

{\subsection{Single-period capacity planning with multiple operational scenarios.}
We consider a single-period planning model for generation and transmission capacity expansion of the Brazilian electricity system targeting the year 2050. All data used to formulate the GenX model is found in~\cite{deng_harmonized_2023} and the related open data repository\footnote{ \url{https://gitlab.com/dlr-ve/esy/open-brazil-energy-data/open-brazilian-energy-data}}. We set a zero-carbon emission cap in 2050, with a penalty for exceeding the cap of $200$ \$/ton CO$_2$. Analogously to the previous section, the penalty for violating multi-day storage state of charge constraints at the beginning/end of any operational sub-period is set equal to twice the cost of non-served energy demand. The dataset~\cite{deng_harmonized_2023} considers $204$ different resource clusters, including 26 hydro-power generators (modeled with multi-day storage constraints), 70 renewable energy generators (i.e. wind and solar), and 81 thermal generators, which we model with linearized unit commitment decisions/constraints. The Brazilian electricity system is subdivided in 27 regions, with 50 lines representing aggregated inter-regional transmission capacity 
Given the level of aggregation in both resource generators and transmission lines, this example considers only continuous investment decisions.

 The dataset~\cite{deng_harmonized_2023} includes hourly time-series for 8 different operational years (2012-2019), based on historical data on renewable energy resources variability, electricity demands, and hydro inflow patterns. As a result, we consider 8 different weather and demand scenarios for 2050, representing uncertainty about future weather conditions. In the absence of a statistical model for these historical time series, we assign the same probability to each of the 8 operational scenarios. Each historical time series for electricity demand is scaled based on the 2050 energy demand projections included in~\cite{deng_harmonized_2023} for scenario labeled as $\texttt{BAU\_V5}$. This process results in $8 \times 52 = 416$ operational sub-periods (e.g. weeks) modeled with hourly resolution, yielding a capacity expansion model with over 53 million variables and 107 million constraints. Note that this model considers less planning decisions compared to the 26-zone multi-period planning case from the previous section, but it includes more operational sub-periods. We utilize this test system to evaluate the performance of regularized Benders decomposition schemes for single-period CEMs with multiple operational scenarios (e.g., stochastic scenarios). We compare the performance of \texttt{Benders-noReg}, $\texttt{Benders-}\Phi^{{\ell_2},(0.2,0.5)}$, $\texttt{Benders-}\Phi^{{int},0.5}$, and $\texttt{Benders-TR}^{0.3}$. In the implementation of the Benders decomposition schemes, we parallelize the solution of the 416 operational sub-problems over $208$ different computing cores, where each core solves two sub-problems in series. 
Table~\ref{tab:results_brazil} reports the result from these computational experiments, where we report the performance only of $\texttt{Benders-}\Phi^{{\ell_2},0.5}$ because it resulted in shorter runtime compared to $\texttt{Benders-}\Phi^{{\ell_2},0.2}$. Again, the interior-point regularized Benders decomposition algorithm $\texttt{Benders-}\Phi^{{int},0.5}$ had the best performance, both in terms of runtime and number of iterations. In particular, $\texttt{Benders-}\Phi^{{int},0.5}$ solved the model in roughly 40 minutes, almost half the runtime required by the basic Benders decomposition algorithm without any regularization. Nonetheless, we observe that \texttt{Benders-noReg} had a better performance compared to Table~\ref{tab:results}, and it even outperformed the $\ell_2-$regularized Benders decomposition algorithm. This indicates that the oscillating behavior typical of Benders decomposition without regularization is more pronounced when the number of planning decisions is larger.
\begin{table}[!h]
    \caption{Results for the Brazilian electricity system. We denote runtime (in minutes) with $t$ and define $k^*$ and $t^*$ as the first iteration and time when the gap goes below 1\%.}
    \label{tab:results_brazil}
    \centering
\begin{tabular}{|l|r|r|r|r|r|}
    \hline
    \text{Method} & \text{Gap (\%)} & \text{Iter} & t (60s) & $k^*$ & $\text{t}^*$ (60s)\\\hline
    $\texttt{Benders-noReg}$ & 0.08 & 55 & 73.04 &   43 & 55.70 \\\hline
    $\texttt{Benders-}\Phi^{{\ell_2},0.5}$ & 0.10 & 75 & 168.22 &   39 & 82.86 \\\hline
    $\texttt{Benders-TR}$ & 0.07 & 52 & 68.43 &   41 & 52.43 \\\hline
    \textbf{$\texttt{Benders-}\Phi^{{int},0.5}$} & \textbf{0.09} & \textbf{33} & \textbf{39.89} &  \textbf{ 25} & \textbf{28.70} \\\hline 
    \end{tabular}
    \end{table}
    
We also solved the full monolithic model with Gurobi, using the barrier method with $\texttt{Method}=2$ and $\texttt{BarConvTol} = 10^{-3}$. The solution of the monolithic model is implemented on a single computing node with 96 cores. Gurobi required 597 minutes to converge (little less than 10 hours). Hence, $\texttt{Benders-}\Phi^{{int},0.5}$ reduced runtime by a factor of 15 compared to this state-of-the-art commercial solver.

Overall, the results reported in this manuscript suggest that the proposed regularized Benders algorithm is effective for solving both multi-period planning models (with foresight) and single-period capacity expansion models with multiple operational scenarios (stochastic operations). Results in Tables~\ref{tab:results} and~\ref{tab:results_brazil} indicate that $\texttt{Benders-}\Phi^{{int},0.5}$ can provide enough computational bandwidth to increase the number of weather and demand scenarios, or further reduce model abstractions and approximation errors, for example increasing spatial resolution, or extending the capacity expansion model to include energy sectors beyond electricity (e.g. natural gas, hydrogen) while maintaining a runtime under 24 hours. 
}

\section{Conclusions}
We have demonstrated that level-set regularization of Benders decomposition methods combined with a formulation suitable to exploit distributed, parallel computation can offer a scalable approach to solve large-scale, mixed-integer capacity expansion models. For a case study with over 111 million variables, our regularized Benders algorithm was able to compute a solution to the continuous relaxation of the planning problem in just over 6 hours and a full solution to the MILP planning problem {(with an optimality gap of $<$0.1\%)} in 8 hours. {We also present a case study of a single-stage planning model with 8 stochastic operational scenarios and solve this continuous problem (with over 53 million variables) in roughly 40 minutes.} Key advantages of the proposed algorithm compared to previously published methods include:
\begin{enumerate}[label=(\roman*)]
    \item Ability to further decompose the operational sub-problem into smaller sub-problems, defined over multi-day sub-periods of each planning stage. This approach can take advantage of parallel computing and improve convergence by adding more cuts per iteration compared to approaches solving year-long operational sub-problems~\cite{Lara2018,li_mixed-integer_2022,zhang_integrated_2023,goke_stabilized_2024}. {This characteristic is particularly well suited to decompose multi-stage planning problems and/or formulations with stochastic operational scenarios, as demonstrated herein.}
    \item Generation of better quality Benders cuts using { planning decisions} in the interior of the feasible set of Problem~\eqref{eq:planning} via level-set regularization, avoiding oscillations and instabilities that affect basic Benders implementations~\cite{jacobson_computationally_2024}.
    \item When solving CEMs with integer planning variables, we initialize Benders iterations by loading pre-computed Benders cuts obtained by solving a continuous relaxation of the problem. For the considered case study, the Benders algorithm with pre-computed cuts starts from a relatively small optimality gap {($\sim$1\%)}, and it converges {to a $<$0.1\% gap} in just a few additional iterations.
\end{enumerate} 
{While including multiple scenarios improves robustness of the planning decisions, the model may still be biased due to the assumptions of perfect foresight and deterministic operational decisions. Note that this limitation affects most open-source capacity expansion models~(e.g.,~\cite{brown_pypsa_2018,johnston_switch_2019,mit_energy_initiative_genx_nodate}). Further research should explore the application of the developed decomposition scheme to stochastic and robust planning formulations accounting for uncertainty in technology costs, policy stringency, or other uncertainties affecting the planning stage decisions. For example, the presented decomposed formulation could be combined with the stochastic dual dynamic programming (SDDP) algorithm~\cite{pereira_multi-stage_1991} to decrease the computational cost of each SDDP iteration. Future work will also investigate the application of our regularized Benders decomposition method to solve capacity expansion models where the operational sub-problems are substituted by a stochastic rolling horizon dispatch procedure as in~\cite{bodal_capacity_2022}.}
Our results show that runtime can still grow rapidly with the size of the system, posing challenges to further increasing spatial/network resolution to accurately represent power flow constraints. Future research should investigate new methods that combine our Benders decomposition framework with network decomposition and partitioning techniques.} 
Finally macro-scale energy systems planning models representing multiple energy networks and industrial supply chains (e.g., TEMOA~\cite{hunter_modeling_2013}, RIO~\cite{williams_carbon-neutral_2021}) present a very similar problem structure as electricity capacity expansion models. The methods presented herein are therefore extensible to this class of models and others with similar problem structure: e.g., investment decisions that must be co-optimized with numerous operational decisions subject to coupling constraints.
\section*{Acknowledgments}
The authors thank Greg Schivley for the generation of the test case data with the package PowerGenome. Funding for this work was provided by the Princeton Carbon Mitigation Initiative (funded by a gift from BP) and the Princeton Zero-carbon Technology Consortium (funded by gifts from GE, Google, ClearPath, and Breakthrough Energy).

\bibliographystyle{IEEEtran}
\bibliography{IEEE_TPS_Benders}


\end{document}